\documentclass[12pt]{article}

\usepackage{amsmath}
\usepackage{latexsym}
\usepackage{amsfonts}
\usepackage{amssymb}

\usepackage{latexsym}
\newcounter{theorem}

\renewcommand{\thetheorem}{\arabic{section}.\arabic{theorem}}

\newenvironment{thm}[1]{\par\addvspace{0.5cm}
    \begin{sloppypar}\refstepcounter{theorem}%
    {\bf #1 \thetheorem.}\it{}}{\end{sloppypar}}

\newcommand{\eh}{\hfill}\newlength{\sperr}

\newenvironment{theorem}{\begin{thm}{Theorem}} {\end{thm}}

\newenvironment{corollary}{\begin{thm}{Corollary}} {\end{thm}}

\newenvironment{defi}[1]{\par\addvspace{0.5cm}
\begin{sloppypar}\refstepcounter{theorem}%
{\bf #1 \thetheorem.}\rm{}}{\end{sloppypar}}

\newenvironment{definition}{\begin{defi}{Definition}}{\end{defi}}

\newenvironment{remark}{\begin{defi}{Remark}}{\end{defi}}

\newcommand{\vs}{\vspace}

\newcommand{\al}{\alpha}
\newcommand{\bt}{\beta}
\newcommand{\dl}{\delta}
\newcommand{\om}{\omega}
\newcommand{\Om}{\Omega}
\newcommand{\lb}{\lambda}

\newcommand{\ve}{\varepsilon}
\newcommand{\gm}{\gamma}
\newcommand{\Gm}{\Gamma}

\newcommand{\vi}{\varphi}

\newcommand{\rn}{\mathbb{R}^n}

\newcommand{\intl}{\int\limits}

\newcommand{\ol}{\overline}

\newcommand{\sbsq}{\subseteq}

\newcommand{\diam}{\mbox{\,\rm diam\,}}

\begin{document}

%\def\theequation{\arabic{section}.\arabic{equation}}
%\renewcommand{\thesection}{\arabic{section}}
%\renewcommand{\theequation}{\arabic{section}.\arabic{equation}}
%\renewcommand{\thetheorem}{\arabic{section}.\arabic{thm}}
%\renewcommand{\thesubsection}{\arabic{section}.\arabic{subsection}}

%-------------------------------------------------------------------------
% editorial commands: to be inserted by the editorial office
%
%\firstpage{1}
%\volume{228}
%\DOI{003-0001}
%\seriesextraline{This is the Concrete Title of this Book\br H.E. R and S.T.C. W, Eds.}
%\firstpage{1}

\centerline{\textbf{Weighted Boundedness of the Maximal,
Singular}}

\centerline{\textbf{and Potential   Operators in   Variable
Exponent Spaces}}

\vspace{4mm} \centerline{by}

\vspace{4mm} \centerline{\textbf{Vakhtang Kokilashvili}}

\centerline{\textbf{Black Sea University,  and  A.Razmadze
Mathematical Institute, Tbilisi, Georgia}}

\centerline{email: \textit{kokil@rmi.acnet.ge}}

\vspace{4mm} \centerline{and}

\centerline{\textbf{Stefan Samko}}

\centerline{\textbf{Universidade do Algarve}}

\centerline{email: \textit{ssamko@ualg.pt}}

\vspace{4mm} AMS Classification 2000: \ \ \ Primary 42B25;
Secondary  47B38

\vspace{4mm} Keywords: maximal functions, potential operators,
weighted Lebesgue spaces, weighted estimates, variable exponent,
metric measure space, doubling condition, Zygmund
 conditions, Zygmund-Bary-Stechkin class

\begin{abstract}
We present a brief survey of recent results on boundedness of some
classical operators within the frameworks of weighted spaces
$L^{p(\cdot)}(\varrho)$ with variable exponent $p(x)$, mainly in
the Euclidean setting and dwell on a new result of the boundedness
of the Hardy-Littlewood maximal operator in the  space
$L^{p(\cdot)}(X,\varrho)$ over a metric measure space $X$
satisfying the doubling condition. In the case where $X$ is
bounded, the weight function satisfies  a certain version of a
general Muckenhoupt-type condition For a bounded or unbounded $X$
we also consider a class of weights of the form
$\varrho(x)=[1+d(x_0,x)]^{\bt_\infty}\prod_{k=1}^m w_k(d(x,x_k))$,
$x_k\in X$, where the functions $w_k(r)$ have finite upper and
lower indices $m(w_k)$ and $M(w_k)$.

Some of the results are new even in the case of constant
$p$.\end{abstract}

\section{Introduction}
\setcounter{equation}{0}\setcounter{theorem}{0}

In this paper, based on the lecture  \textit{Harmonic Analysis on
variable exponent spaces} presented by the second author at AMADE
International Conference AMADE-2006 (Analytic Methods of Analysis
and Differential Equations) held in   "Staiki", Minsk, \  Belarus,
September 13-19, 2006, we present a survey of a certain selection
of known results  on weighted estimations of operators in variable
exponent Lebesgue spaces  and expose some new results.

The study of classical operators of harmonic analysis (maximal,
singular operators and potential type operators) in the
generalized Lebesgue spaces $L^{p(\cdot)}$ with variable exponent,
weighted or non-weighted, undertaken last decade, continues to
attract a strong interest of researchers, influenced in particular
by possible applications revealed in the book \cite{525}.
 We refer in particular to the surveying articles \cite{106b}, \cite{316b}, \cite{580bd}.

The area which is now called variable exponent analysis, last
decade became a rather branched field with many interesting
results obtained in Harmonic Analysis, Approximation Theory,
Operator Theory, Pseudo-Differential Operators, the survey of
which would be a good task. Our survey is far from being complete
in this sense. In our selection of papers for this survey  we have
mainly chosen results on weighted estimates, obtained  after the
above surveys had appeared,  and present some  new results on
weighted boundedness of the Hardy-Littlewood maximal
 operator in such spaces over metric measure spaces.

 \vspace{4mm} The paper is organized as follows.  In Section \ref{tri} we recall some basics from the
theory of generalized Lebesgue spaces with variable exponent on
metric measure spaces. In Section \ref{Phi} we give necessary
preliminaries  on  the upper and lower indices (of
Matuszewska-Orlicz-type) of weights in the Zygmund-Bary-Stechkin
class.  In section \ref{dobavtri} we give a survey of recent
results  on the weighted boundedness of classical operators of
harmonic analysis (Hardy-Littlewood maximal operator,
Calder\'on-Zygmund type singular operators, Cauchy singular
integral operator on Carleson curves, potential type operators and
some classes of convolution operators) in weighted Lebesgue spaces
with variable exponent.

 In Section \ref{yyyys4} we give new
results - Theorems A, B and C - on the weighted boundedness of the
maximal operator on metric measure spaces with doubling condition.

Theorem A gives a kind of "Muckenhoupt-looking" condition. In
Theorem B we deal with a bounded space $X$ and radial-type
oscillating weights $w[d(x_0,x)]$; we obtain sufficient conditions
for the boundedness in this case in the form
\begin{equation}\label{avotyavamhu}
-\frac{m(\mu B)}{p(x_0)} < m(w)\le M(w) < \frac{m(\mu
B)}{p^\prime(x_0)},
\end{equation}
where \begin{equation}\label{introduced}
 m(\mu B)= \lim_{t\to 0}\frac{\ln\left(\limsup\limits_{r\to
0}\inf\limits_{x\in X} \frac{\mu B(x,rt)}{\mu B(x,r)}\right)}{\ln t}
\end{equation}
 is a kind of
"uniform" lower Matuszewska-Orlicz-type index of the function $\mu B(x,r)$, $ \ B(x,r)=\{y\in X:
d(x,y)< r\}$ with respect to the variable $r$: see Section \ref{Phi}  for Matuszewska-Orlicz-type
indices where one may find hints onto how these indices of the measure $\mu B(x,r)$ appeared. In
Theorem C we give a version of Theorem B for the case of unbounded spaces $X$.

\vspace{6mm}

\noindent \textbf{Notation}

\vspace{3mm}
 \noindent
$(X,d,\mu)$ is a  measure space with quasimetric $d$  and a
non-negative measure $\mu$;

\noindent $B(x,r)=B_X(x,r)=\{y\in X: d(x,y)<r\}$;

\noindent $p^\prime(x)=\frac{p(x)}{p(x)-1}$, \ \ $1 <
p(x)<\infty$, \ \ $\frac{1}{p(x)}+ \frac{1}{p^\prime(x)}\equiv 1$
;

\noindent $p_-=p_-(X)=\inf\limits_{x\in X} p(x)$,
$p^+=p^+(X)=\sup\limits_{x\in X} p(x)$;

\noindent $p^\prime_-=\inf\limits_{x\in X}
p^\prime(x)=\frac{p^+}{p^+-1}$, \ $(p^\prime)^+=\sup\limits_{x\in
X} p^\prime(x)=\frac{p_-}{p_--1}$;

\noindent $\mathbb{P}(X)$, see (\ref{1st})-(\ref{2nd});

\noindent $C,c$ may denote different positive constants;

 \noindent a.i. =almost increasing
$\Longleftrightarrow u(x)\le Cu(y)$ for $x\le y, C>0$.

\vspace{3mm}

\section{Some basics for variable exponent spaces}\label{tri}

\setcounter{theorem}{0} \setcounter{equation}{0}

In the sequel $(X,d,\mu)$ is a homogeneous type space, i.e. a
measure space with quasimetric $d$  and a non-negative measure
$\mu$ satisfying the doubling condition; we refer to \cite{187},
\cite{221a}, \cite{225a} for the basic notions of function spaces
on metric measure spaces. We suppose that the measure $\mu$
satisfies the doubling condition
\begin{equation}\label{double}
\mu B(x,2r) \le C \mu B(x,r).
\end{equation}

  By $\mathbb{P}(X)$ we denote the set of bounded measurable functions  $p(x)$ defined on $X$ which
   satisfy
the conditions
\begin{equation}\label{1st}
1<p_-\le p(x)\le p^+<\infty, \;\; x\in X
\end{equation}
and
\begin{equation}\label{2nd}
|p(x)-p(y)|\le \frac{A}{\ln\frac{1}{d(x,y)}}\,, \;\; d(x,y)\le
\frac{1}{2}, \;\; x,\,y\in X.
\end{equation}

By $L^{p(\cdot)}(X,\varrho)$, where $\varrho(x)\ge 0$, we denote
the weighted Banach space of measurable functions $f: X\to
\mathbb{C}$ such that
\begin{equation}\label{d4cd}
\|f\|_{L^{p(\cdot)}(X,\varrho)}:=\|\varrho
f\|_{p(\cdot)}=\inf\left\{\lb>0: \intl_X
\left|\frac{\varrho(x)f(x)}{\lb}\right|^{p(x)} \;d\mu(x)\le
1\right\}<\infty .
\end{equation}
We write $L^{p(\cdot)}(X,1)=L^{p(\cdot)}(X)$  and $\|f\|_{L^{p(\cdot)}(X)}=\|f\|_{p(\cdot)}$ in the
case $\varrho(x)\equiv 1$.

The generalized Lebesgue  spaces $L^{p(\cdot)}(X)$ and Sobolev
spaces $W^{1,p(\cdot)}$ with variable exponent on metric measure
spaces have been considered in \cite{145c}, \cite{167z},
\cite{167b}, \cite{224ab}, \cite{224b}, \cite{224ba}, \cite{307b},
\cite{405a}, the Euclidean case being studied in \cite{146},
\cite{160a}, \cite{332}, \cite{618a}, see also references therein.
     We recall  the H\"{o}lder inequality
%   (3.2)
\begin{equation}
\int\limits_X |f(x)g(x)|\,d\mu(x)\le
k\big\|f\big\|_{p(\cdot)}\cdot \big\|g \big\|_{p^\prime(\cdot)}
\end{equation}
where $k=\frac{1}{p_-}+\frac{1}{p^\prime_-}$. We note also that
the embedding
\begin{equation}\label{35}
 L^{p(\cdot)}\sbsq
L^{s(\cdot)}, \;\; \; \big\|f\big\|_{s(\cdot)}\le
C\big\|f\big\|_{p(\cdot)},
\end{equation}
 is valid for $1\le s(x)\le p(x)\le
p^+<\infty$, when $\mu (X) <\infty$.

\section{Preliminaries on Zygmund-Bary-Stechkin classes.}\label{Phi}
\setcounter{theorem}{0} \setcounter{equation}{0}

In this section we follow some ideas of papers \cite{270a},
\cite{539}, \cite{539d}, \cite{539h} . Let $0<\ell \le \infty$.
 We denote
\begin{equation}\label{f:2.1abc} W=\{w\in C([0,\ell]):  \
w(t)>0 \ \ \textrm{for}\ \ \ t>0, \ \ w(t) \ \ \textrm{is almost
increasing}\}
\end{equation}
and
\begin{equation}\label{f:2.1abcd}
W_0=\{w\in W: w(0)=0\}.
\end{equation}
 We also need a wider class
\begin{equation}\label{55110m}
\widetilde{W}([0,\ell])=\{\varphi: \ \ \exists
a=a(\vi)\in\mathbb{R}^1 \ \ \textrm{such that}\ \ \ x^a\vi(x)\in
W([0,\ell])\}.
\end{equation}

\subsection{The Zygmund-Bary-Stechkin type classes
$\Phi^\al_\bt=\Phi^\al_\bt([0,\ell])$ and
$\Psi^\al_\bt=\Psi^\al_\bt([\ell,\infty]),\  0<\ell<\infty$.}
\label{hnb7hx4}

 The following class $\Phi_\bt^\al$ of Zygmund-Bary-Stechkin type in the case $\al=0$ and $ \bt=1,2,3,...$ was introduced in
  \cite{46} (in \cite{46} functions   $w$ were assumed to be  increasing functions).
\begin{definition}\label{def:sew2kdd} The  Zygmund-Bary-Stechkin type class $\Phi_\bt^\al=\Phi_\bt^\al([0,\ell]), \ -\infty< \al<\bt <\infty,$
is defined  as $\Phi_\bt^\al:= \mathcal{Z}^\al\cap
\mathcal{Z}_\bt$, where $\mathcal{Z}^\al$ is the class of
functions $w\in \widetilde{W}$ satisfying  the
 condition
$$
 \int_0^h\frac{w(t)}{t^{1+\al}}dt
\le c\frac{w(h)}{h^\al} \eqno(\mathbb{Z}^\al)
$$
and $\mathcal{Z}_\bt$ is the class of functions $w\in W$
satisfying the
 condition
$$
 \int_h^\ell\frac{w(t)}{t^{1+\bt}}d(t) \le c\frac{w(h)}{h^\bt},
 \eqno(\mathbb{Z}_\bt)$$
where $c=c(w)>0$ does not depend on $h\in (0,\ell].$
\end{definition}

\vspace{5mm}

We also need a  class of functions with a similar behaviour at
infinity. Let
 $C_+([\ell,\infty)), \ 0<\ell<\infty$, be  the
class of functions $w(t)$ on $[\ell,\infty)$, continuous and
positive at every point $t\in[\ell,\infty)$ \ and having a finite
or infinite limit $\lim\limits_{t\to \infty}w(t)=: w(\infty)$. We
denote
\begin{equation}\label{f:2.1abchiuw}
W=W([\ell,\infty))=\{w\in C_+([\ell,\infty):  \
 w(t) \ \ \textrm{is a.i.}\}
\end{equation}
and
\begin{equation}\label{55110mnewa}
\widetilde{W}([\ell,\infty))=\{\varphi: \ \ \exists
a=a(\vi)\in\mathbb{R}^1 \ \ \textrm{such that}\ \ \ x^a\vi(x)\in
W([\ell,\infty)\}.
\end{equation}

\begin{definition}\label{def:sew2kddsa3} Let $ \ -\infty< \al<\bt
<\infty.$ We put   $\Psi^\bt_\al:= \widehat{\mathcal{Z}}^\bt\cap
\widehat{\mathcal{Z}}_\al,$ where $\widehat{\mathcal{Z}}^\bt$ is
the class of functions $w\in \widetilde{W}([\ell,\infty))$
satisfying  the
 condition
\begin{equation}\label{asw8n}
 \int_r^\infty\left(\frac{r}{t}\right)^\bt \frac{w(t)\, dt}{t}
\le c w(r), \ \ \ \  r \to \infty,
\end{equation}and $\widehat{\mathcal{Z}}_\al$ is the class of functions $w\in
W([\ell,\infty))$ satisfying the
 condition
\begin{equation}\label{asw9n}
 \int_\ell^r\left(\frac{r}{t}\right)^\al\frac{w(t)\, dt}{t} \le
 cw(r),\ \ \ \ \ r \to \infty
\end{equation}
where $c=c(w)>0$ does not depend on $r\in [\ell,
\infty).$
\end{definition}

Observe that properties of functions in the class
$\Psi^{\bt}_{\al} ([\ell,\infty))$ are easily derived from those
of functions in $\Phi^{\al}_{\bt} ([0,\ell])$ because of the
following equivalence
\begin{equation} \label{ggs99}
w\in \Psi^{\bt}_{\al} ([\ell,\infty)) \  \ \ \Longleftrightarrow \
\ \ w_\ast \in \Phi^{-\bt}_{-\al}([0,\ell^\ast]),
\end{equation}
where $w_\ast(t)=w\left(\frac{1}{t}\right)$ and
$\ell_\ast=\frac{1}{\ell}.$

\subsection{Index numbers $m(w)$ and $M(w)$ of non-negative a. i. functions}\label{hhhx4}

The numbers
\begin{equation}\label{m}
m(w)=\sup_{t>1}\frac{\ln\left(\liminf\limits_{h\to 0}
\frac{w(ht)}{w(h)}\right)}{\ln
t}=\sup_{0<t<1}\frac{\ln\left(\limsup\limits_{h\to 0}
\frac{w(ht)}{w(h)}\right)}{\ln t}= \lim_{t\to
0}\frac{\ln\left(\limsup\limits_{h\to 0}
\frac{w(ht)}{w(h)}\right)}{\ln t}
\end{equation}
 and
\begin{equation}\label{M}
M(w)=\inf_{t>1}\frac{\ln\left(\limsup\limits_{h\to 0}
\frac{w(ht)}{w(h)}\right)}{\ln t}=\lim_{t\to
\infty}\frac{\ln\left(\limsup\limits_{h\to 0}
\frac{w(ht)}{w(h)}\right)}{\ln t}
\end{equation}
 (see \cite{539}, \cite{539e}), \cite{539d}, will be referred to  as \textit{the lower and upper indices} of
the function $w(t)$. We have $0\le m(w)\le M(w) \le \infty  \ \
\textrm{for}\ \  w\in W$.

\vspace{2mm} The  indices $m(\om)$ and $M(\om)$ may be also well
defined for functions $w(t)$ positive for $t>0$ which do not
necessarily belong to $W$, for example, for $w\in \widetilde{W}$.
Observe that
\begin{equation}\label{men}
m(w_a)=a+m(w), \ \ \ \ M(m_a)=a+M(w) \ \ \ \ \textrm{where}      \
\ \ \ \ w_a(t):=t^aw(t)
\end{equation}
and
\begin{equation}\label{menajsl}
m(w^\lb)=\lb m(w), \quad M(w^\lb)=\lb M(w), \quad \lb\ge 0
\end{equation}
for every $w\in\widetilde{W}$.

\vspace{4mm}

The indices $m_\infty(w)$ and $M_\infty(w)$ responsible for the behavior of functions $w\in
\Psi^{\bt}_{\al} ([\ell,\infty))$ at infinity are introduced in the way similar to (\ref{m}) and
(\ref{M}):
\begin{equation}\label{msusut}
 m_\infty(w) =\sup_{x>1}\frac{\ln \
\left[\liminf\limits_{h\to \infty} \frac{w(xh)}{w(h)} \right]}{\ln
\ x}\ , \ \  M_\infty(w) =\inf_{x>1}\frac{\ln \
\left[\limsup\limits_{h\to \infty} \frac{w(xh)}{w(h)} \right]}{\ln
\ x}.
\end{equation}

\section{A survey of recent results on  boundedness of classical operators in weighted spaces
$L^{p(\cdot)}(\Om,\varrho)$ }\label{dobavtri}

\setcounter{theorem}{0} \setcounter{equation}{0}

In this section we consider the following classical operators: \\

\noindent \textit{1) \ Convolution operators}
$$Af(x)=\intl_{\mathbb{R}^n}k(y)f(x-y)dy$$
with rather "nice" kernels without log-condition,

\noindent \textit{2) \ Hardy-Littlewood maximal operator}
\begin{equation}\label{weightmax}
\mathcal{M} f(x)=\sup_{r>0} \frac{1}{\mu (B(x,r))}
\int\limits_{B(x,r)} |f(y)|\,d\mu(y), \quad x\in X
\end{equation}
where $X$ in general is a metric measure space, being either an
open set in  $\rn$ or a Carleson curve on the complex plane in
this section.

\noindent \textit{3) \ Calder\'on-Zygmund type singular integral
operator}
\begin{equation}\label{f:5sd}
T f(x)=\lim\limits_{\varepsilon\to
0}\int\limits_{|x-y|>\varepsilon} k(x,y) f(y)\; dy
\end{equation}
with  the so called \textit{standard kernel} (see, for instance,
\cite{132}, p.99), and also the Cauchy singular integral
\begin{equation}\label{f2axsdf20}
S_\Gm f(t)=\frac{1}{\pi i}\intl_{\Gm}
\frac{f(\tau)}{\tau-t}d\nu(\tau)
\end{equation}
along Carleson curves $\Gm$ on complex plane, where $\nu$ is the
arc-length measure; we recall that $\Gm$ is called a
\textit{Carleson curve}, if it satisfies the condition
$$\nu(\Gm\cap B(t,r))\le Cr$$
where the constant $C>0$ does not depend on $t\in\Gm$ and $r>0$;

\noindent \textit{4) \ potential type operators}
\begin{equation} \label{potentialw}
I^{\al(\cdot)} f(x) \ = \  \int\limits_{\Om}\frac{f(y)\;
dy}{|x-y|^{n-\al(x)}} , \ \ \ \ \ \ 0<\inf\al(x)\le \sup \al(x)< n
\end{equation}
over open bounded sets $\Om$ in $\rn$, and

\noindent \textit{5) \ Hardy operators}
\begin{equation} \label{0v08ccc9}
H^\al f(x)= x^{\al-1}\intl_{0}^x\frac{f(y)}{y^\al}dy \ \ \ \ \ \
\textrm{and}\ \ \ \ \ \ \ \mathcal{H}_\bt f(x)=
x^{\bt}\intl_x^\infty\frac{\vi(y)\,dy}{y^{\bt +1}}.
\end{equation}

Observe that boundedness of various classical operators in the
non-weighted case was proved in \cite{101zb} by the extrapolation
method.  In relation to the extrapolation method we refer also to
\cite{101ac},  \cite{101b}, \cite{101ad}.

\subsection{On convolution operators}\label{convoluton}

As is well known, the Young theorem in its natural form is not
valid in the case of variable exponent, whatsoever smooth exponent
$p(x)$ is. However, a natural expectation was that the Young
theorem may be valid in the case of rather "nice" kernels and even
without log-condition. This is true indeed, see Theorem \ref{5556}
below.

By $\mathcal{P}_{\infty}(\mathbb{R}^n)$
 we denote the set  of measurable bounded functions $p: \mathbb{R}^n\to \mathbb{R}^1_+$
which satisfy the following conditions:\\
$ i) \ \ 1\le p_-\le p(x)\le p_+<\infty, \ \ \ \ x\in
\mathbb{R}^n$,\\
 $ ii)$ \ there exists $p(\infty)=\lim\limits_{x\to \infty}p(x)$ \ \ \
and
\begin{equation}\label{55sa7}
|p(x)-p(\infty)|\le \frac{A}{\ln\,(2+|x|)}, \ \ x\in \mathbb{R}^n.
\end{equation}

The following statement was proved in \cite{107da}.
\begin{theorem}\label{5556}
Let $k(y)$ satisfy the estimate
\begin{equation}\label{eq:decrease}
|k(y)|\le \frac{C}{(1+|y|)^\lb}, \ \ \ \ \ \ \ y\in \mathbb{R}^n
\end{equation}
for some $\lb >
n\left(1-\frac{1}{p(\infty)}+\frac{1}{q(\infty)}\right).$ Then the
convolution operator
 is bounded from the
space $L^{p(\cdot)}(\mathbb{R}^n)$ to the space
$L^{q(\cdot)}(\mathbb{R}^n)$ under the only assumption that
$p,q\in \mathcal{P}_\infty(\mathbb{R}^n)$ and $q(\infty)\ge
p(\infty)$.
\end{theorem}

The convergence in $L^{p(\cdot)}(\mathbb{R}^n)$ of convolutions with identity approximation kernels
was studied in \cite{101zza}.

\subsection{On Maximal Operator}\label{subsmaxop}

 Non-weighted boundedness of the maximal operator was first proved in  \cite{104b}, \cite{106}
 for bounded domains or for $\mathbb{R}^n$ with $p(x)\equiv const$
 outside some large ball.  For further results in non-weighted case see
\cite{101ab}, \cite{101aba}, \cite{106z}, \cite{355a}, \cite{414b}, \cite{414d}. A special
situation when $p(x)$ may tend to 1 or $n$, was studied in \cite{101zz},  \cite{167a}, \cite{
224zb}, \cite{224c}. For the non-validity of the modular inequality for the maximal function in
case of non-constant $p(x)$ we refer to \cite{355b}.

The result of \cite{104b}, \cite{106}  was extended to the weighted case with power weights in
\cite{321a}. After that several papers were devoted to consideration of more general weights for
the variable exponent setting.

 A characterization of general weights admissible for the boundedness of the maximal operator in the spirit of
  Muckenhoupt type condition is still an open question.  The generalization to the case
of more general weights encountered essential difficulties.

 The main progress  in obtaining sufficient conditions on weights, similar to those for constant $p$,
  in the Euclidean setting was obtained for a certain special class of
 weights, although essentially more general than the class of  power weights.
 This class consists of  weights of the form
\begin{equation}\label{ves}
\varrho(x) = [1+w(|x|)]\prod\limits_{k=1}^m w_k(|x-x_k|), \ \ \ \
x_k\in\overline{\Om} \subset \rn
\end{equation}
    of radial-type, the factor $1+w(|x|)$  appearing in the case of unbounded
sets $\Om$. These weights have a typical feature of Muckenhoupt
weights: they may oscillate between two power functions. This
class may be also interpreted as a kind of Zygmund-Bary-Stechkin
class. The introduction of this class of weights is based on the
observation that the integral constructions involved in the
Muckenhoupt condition  for radial weights (in the case of constant
$p$) are exactly those which appear in the Zygmund-type
conditions.

 For  weights in this class it proved to be possible to obtain  sufficient conditions of the
 boundedness of the maximal operator in terms of the so called upper and lower index
numbers $m(w_k)$ and $M\left(w_k\right)$ of the weights $w_k(r)$ (Matuszewska-Orlicz type indices).
These conditions are obtained in the form of the  natural numerical intervals
\begin{equation}\label{prov}
 -\frac{n}{p(x_k)} <m(w_k)\le M(w_k)<\frac{n}{p^\prime(x_k)}
\end{equation}
   "localized" to the nodes $x_k$ of
the weights $w_k(|x-x_0|)$. The sufficiency of this condition in
terms of the numbers $m(w)$ and $M\left(w\right)$ is a new result
even in the case of constant $p$. As is known, even in the case of
constant $p$  the verification of the Muckenhoupt condition for a
concrete weight may be an uneasy task.  Therefore, independently
of finding an analogue of the Muckenhoupt condition for variable
exponents, it is always of importance to find easy to check
sufficient conditions for weight functions, as for instance in
(\ref{prov}). The following theorem was proved in \cite{317c}, see
also a sketch of the proof in \cite{321i}. The
Zygmund-Bary-Stechkin class $\Phi_n^0$ and the notion of index
numbers $m(w)$ and $M(w)$ are defined in Section \ref{Phi}.
\begin{theorem}\label{theorem}
 Let $X=\Om$ be a bounded domain in
$\mathbb{R}^n$ with Lebesgue measure,  let  $p\in \mathbb{P}(\Om)$
and $\varrho$ be weight (\ref{ves}) with $w \equiv 0$. The
operator $\mathcal{M}$ is bounded in $L^{p(\cdot)}(\Om,\varrho)$,
if $r^{\frac{n}{p(x_k)}}w_k(r)\in \Phi_n^0$ or, which is
equivalent, conditions (\ref{prov}) are satisfied.
\end{theorem}

\vs{3mm} In the case of power weights, a similar statement was proved in the context of  Carleson
curves $\Gm$ on complex planes, which are examples of metric measure spaces with arc length measure
and coinciding upper and lower dimensions equal to $1$. Because of the interest to the case $X=\Gm$
in the theory of singular integral equations, we formulate separately this statement proved in
\cite{321j} (see also a sketch of the proof in \cite{321i}), including  the case of infinite curves
$\Gm$.

In case $X=\Gm$, for points on $\Gm$ we agree to write $t$ instead
of $x$, $t_k$ instead of $x_k$, etc

\begin{theorem}\label{theorem1}  Let
 \ $\Gm$ be a simple Carleson curve of finite or infinite
 length, let
  $p\in \mathbb{P}(\Gm)$ and  $p(t)\equiv
p_\infty =const$ for $t\in \Gm\backslash (\Gm\cap B(0,R))$ for
some $R>0$, if \ $\Gm$  is infinite.
 Then the maximal operator $\mathcal{M}$
is bounded in the space $L^{p(\cdot)}(\Gm,\varrho)$ with weight
\begin{equation}\label{vespower}
\varrho(t) = (1+|t|)^\bt\prod\limits_{k=1}^m |t-t_k|^{\bt_k}, \
t_k\in\Gm,
\end{equation}
 if and only if
\begin{equation}\label{condit}
 -\frac{1}{p(t_k)}<\bt_k<\frac{1}{p^\prime(t_k)}, \ \ \  \ \
k=1,...,m, \ \ \ \  \textrm{and} \ \ \ \
-\frac{1}{p_\infty}<\bt+\sum\limits_{k=1}^m\bt_k
<\frac{1}{p^\prime_\infty},
\end{equation}
 the latter condition appearing in the
case $\Gm$ is infinite.
\end{theorem}

We conclude this subsection by the observation  that the
fractional maximal operator for variable exponent setting was
studied in \cite{ 73a}, \cite{101zb} and \cite{321ea}.

\subsection{Weighted estimates for Singular Operators}\label{subssingul}

\vs{4mm} For non-weighted estimates of singular integrals we refer to \cite{107a}.

\vs{4mm}\textbf{a) \ Calder\'on-Zygmund type operators.} The
following theorem on weighted boundedness of Calder\'on-Zygmund
type singular operators (\ref{f:5sd}) was proved in \cite{317e}.

\begin{theorem}\label{tehw3d}
 Let $\Om$ be a bounded open set in $\rn$ and $p\in \mathbb{P}(\Om)$.
A singular operator $T_\Om$ with a standard kernel $k(x,y)$,
bounded from $L^1(\Om)$ to $L^{1,\infty}(\Om)$, is bounded in the
space $L^{p(\cdot)}(\Om,\varrho)$ with the weight $
\varrho(x)=\prod_{k=1}^mw_k(|x-x_k|), $
 where $x_k\in \ol{\Om}$ and
\begin{equation}\label{enrey}
w_k(r), \frac{1}{w_k(r)}\in \widetilde{W}([0,\ell]), \ \ \
k=1,..., m, \quad \textrm{and} \quad \ell=\diam \Om,
\end{equation}
if conditions (\ref{prov}) are satisfied.
\end{theorem}

Theorem \ref{tehw3d} was obtained in \cite{317e} by means of  the
following general result, where
\begin{equation}\label{muckenh}
 \mathcal{A}_{p(\cdot)}(\rn) = \{\varrho : \
 \textrm{the maximal operator}\  \mathcal{M} \quad \textrm{is
bounded in} \quad  L^{p(\cdot)}(\rn,\varrho)\}
\end{equation}

\begin{theorem}\label{7} (\cite{317e})
Let $p\in \mathbb{P}(\rn)$and let the weight function
$\varrho$ satisfy the assumptions\\
$i) \ \ \varrho \in \mathcal{A}_{p(\cdot)}(\rn) \ \ \ \textrm{and}
\ \ \
\frac{1}{\varrho}\in\mathcal{A}_{p^\prime(\cdot)}(\mathbb{R}^n)$;\\
$ii) \ \ \textrm{there exists an}  \ s\in (0,1)  \ \ \textrm{such
that} \ \
\frac{1}{\varrho^s}\in\mathcal{A}_{\left({\frac{p(\cdot)}{s}}\right)^\prime}(\mathbb{R}^n)$.
\\  Then a singular operator $T$ with a standard kernel $k(x,y)$ and
of weak (1,1)-type, is bounded in the space
$L^{p(\cdot)}(\mathbb{R}^n,\varrho)$.
\end{theorem}

The non-weighted case $\varrho\equiv 1$ of Theorem \ref{7}  was
proved in \cite{107a}.

 \vs{4mm} The following theorem
was proved in \cite{503a} in the context adjusted to the theory of
PDO for operators of the form
\begin{equation}\label{f3v76a}
\mathbb{A}f(x)=\intl_{\mathbb{R}^n}k(x,x-y)f(y)dy
\end{equation}
where $k(x,z)\in C^1(\mathbb{R}^n\times (\mathbb{R}^n\backslash
\{0\}))$ and it is assumed that
\begin{equation}\label{f2ava}
\lb_1(\mathbb{A}):=\sup\limits_{|\al|=1}\sup\limits_{x,z\in\mathbb{R}^n\times
\mathbb{R}^n}|z|^{n+1}\left|\partial_x^\al k(x,z)\right| <\infty,
\end{equation}
\begin{equation}\label{f2avaxt}
 \lb_2(\mathbb{A}):=\sup\limits_{|\bt|=1}\sup\limits_{x,z\in\mathbb{R}^n\times
\mathbb{R}^n}|z|^{n+1}\left|\partial_z^\bt k(x,z)\right| <\infty
\end{equation}
and the operator $\mathbb{A}$ is of weak (1,1)-type:
\begin{equation}\label{5ossss}
|\{x\in \mathbb{R}^n: |\mathbb{A}f(x)|>t\}|\le
\frac{C(\mathbb{A})}{t}\intl_{\mathbb{R}^n}|f(x)|\,dx.
\end{equation}

 We consider  power weights
\begin{equation}\label{vespowerrn}
\varrho(x) = (1+|x|)^\bt\prod\limits_{k=1}^m |x-x_k|^{\bt_k}, \
x_k\in\rn.
\end{equation}

\begin{theorem} \label{fffa}
 Let the operator
$\mathbb{A}$ satisfy conditions (\ref{f2ava})-(\ref{5ossss}). \\
I. Let \ $p\in\mathbb{P}(\rn)$ satisfy the decay condition
\begin{equation}\label{fcvc4p}
|p(x)-p(\infty)|\le \frac{A}{\ln\;(2+|x|)}, \ \ \ \ \ \
x\in\mathbb{R}^n,
\end{equation}
 Then the operator $\mathbb{A}$ is bounded in the space
$L^{p(\cdot)}(\mathbb{R}^n)$.\\
II.  Let \  $p\in\mathbb{P}(\rn)$ be constant at infinity:
$p(x)\equiv const=p_\infty$ for $|x|\ge R$ with some $R>0$.
 Then the operator  $\mathbb{A}$ is bounded in the space
$L^{p(\cdot)}(\mathbb{R}^n,\varrho)$  with weight
(\ref{vespowerrn}), if
\begin{equation}\label{f3}
-\frac{n}{p(x_k)}<\bt_k <  \frac{n}{p^\prime(x_k)}, \ \ \
k=1,...,m, \quad \textrm{and}
-\frac{n}{p_\infty}<\bt+\sum\limits_{k=1}^m\bt_k <
\frac{n}{p^\prime_\infty}.
\end{equation}
In both cases I and II
\begin{equation}\label{fc888s}
\|\mathbb{A}\|_{L^{p(\cdot)}(\mathbb{R}^n,\varrho)}\le
c(n,p,\varrho)\left[\lb_1(\mathbb{A})+\lb_2(\mathbb{A})+C(\mathbb{A})\right]
\end{equation}
where the constant  $c(n,p,\varrho)$ depends only on $n, p$ and
$\varrho$.
\end{theorem}

For pseudo-differential operators
\begin{equation}\label{e1.2}
Au(x)=\left( 2\pi \right) ^{-n}\int_{\mathbb{R}^{n}}d\xi \int_{\mathbb{R}%
^{n}}a(x,\xi )u(y)e^{i(x-y,\xi )}dy
\end{equation}
we obtain the following corollary in which the known L.
H\"{o}rmander class $ S_{1,0}^{0}$ is the class of  functions
 $a\in C^{\infty }\left( \mathbb{R}_{x}^{n}\times \mathbb{R}_{\xi }^{n}\right) ,$
 such that
\begin{equation}
\left\vert a\right\vert _{r,t}=\sum_{\left\vert \alpha \right\vert
\leq
r,\left\vert \beta \right\vert \leq t}\sup_{\mathbb{R}^{n}\times \mathbb{R}%
^{n}}\left\vert \partial _{\xi }^{\alpha }\partial _{x}^{\beta
}a(x,\xi )\right\vert \left\langle \xi \right\rangle ^{\left\vert
\alpha \right\vert }<\infty  \label{1.1}
\end{equation}%
for all the multi-indices $\alpha ,\beta $.

\begin{corollary}\label{cor3}
Let  $p\in\mathbb{P}(\rn)$ satisfy the decay condition
(\ref{fcvc4p}). Then pseudo differential  operators $A$ with
symbols $a(x,\xi )$ in the class $ S_{1,0}^{0}$ are bounded in the
space $L^{p(\cdot )}(\mathbb{R}^{n}).$
\end{corollary}

These results were used in \cite{503a} to establish criteria for
Fredhomness of pseudo differential operators.

\vs{4mm}\textbf{b) \ Cauchy singular integral operator.}
 A theorem
on the boundedness  of the Cauchy singular operator in the
variable exponent setting was first obtained for rather smooth
curves, namely, Lyapunov curves or curves of bounded turning
without cusps, see \cite{321c}. Meanwhile the modern development
of the operator theory related to singular integral equations
required  a validity of such a boundedness on an arbitrary
Carleson curve.
 The following result was proved in
\cite{317b}.

\begin{theorem}\label{theorem2} Let
 \ $\Gm$ be a simple Carleson curve of  finite or infinite
 length, let
  $p\in \mathbb{P}(\Gm)$
 and the
following condition at infinity
$$
|p(t)-p(\tau)|\le
\frac{A_\infty}{\ln\frac{1}{\left|\frac{1}{t}-\frac{1}{\tau}\right|}},
\ \ \ \ \ \ \left|\frac{1}{t}-\frac{1}{\tau}\right|\le
\frac{1}{2},
$$
for $ |t|\ge L, \ |\tau|\ge L$ with some $L>0$, in the case $\Gm$
is an infinite curve.
 Then the singular operator $S_\Gm$ is bounded in the space
$L^{p(\cdot)}(\Gm,\varrho)$  with weight (\ref{vespower}), if and
only if conditions (\ref{condit}) are satisfied.
\end{theorem}

An extension of Theorem \ref{theorem2} to the case of oscillating
weights from the Zygmund-Bary-Stechkin class $\Phi^\bt_\dl$ was
given in \cite{317e}.

We also mention the following boundedness result, obtained in
\cite{317ab}, admitting a wide class of oscillating weights; it is
an extension to the case of variable exponents of the known
Helson-Szeg\"o theorem \cite{226a}. To this end we need the
notation
$$W^{p(\cdot)}(\Gm)=\left\{\varrho: \ \varrho S_\Gm \frac{1}{\varrho} \ \ \  \textrm{is bounded in} \ \ \
L^{p(\cdot)}(\Gm)\right\}.$$

 \begin{theorem}\label{c6.1} Let $\Gm$ be a bounded
Carleson curve and  $p\in \mathbb{P}$ , $ \varrho\in
W^{p(\cdot)}(\Gm)$ and $\frac{1}{\varrho}\in
L^{p^\prime(\cdot)+\ve}(\Gamma)$, where $\ve
>0$. \  Then  the function
$$\varrho_\vi(t)=\varrho(t)\exp{\left|\frac{1}{2\pi}\intl_{\Gm}\frac{\vi(\tau)\; d\tau}{\tau-t}\right|}$$
 with real continuous $\vi$, belongs to $W^{p(\cdot)}(\Gm)$.
\end{theorem}

This statement, as formulated  in Theorem \ref{c6.1}, follows from
Corollary 6.2 in \cite{317ab}, if we take into account that for a
variable exponent $p\in \mathbb{P}$ the class of curves $\Gm$ for
which the singular operator $S_\Gm$ is bounded in
$L^{p(\cdot)}(\Gm)$, coincides with the class of Carleson curves,
as shows the following theorem, proved  in \cite{317b}.

\begin{theorem} Let $\Gm$ be a finite rectifiable curve.
Let $p: \Gm\to [1,\infty)$  be a continuous function. If the
singular operator $S_\Gm$ is bounded in the space
$L^{p(\cdot)}(\Gm)$, then the curve $\Gm$ has the property
\begin{equation}\label{f3nuskus2}
\sup\limits_{t\in\Gm\atop r>0}\frac{\nu(\Gm\cap
B(t,r))}{r^{1-\ve}}<\infty
\end{equation}
for every $\ve>0$. If $p(t)$ satisfies  the log-condition
(\ref{2nd}), then property (\ref{f3nuskus2}) holds with $\ve=0$,
that is, $\Gm$ is a Carleson curve.
\end{theorem}

\vs{4mm} From Theorem \ref{theorem2} the following statement
important for applications is derived, see \cite{317b}.
\begin{theorem}\label{theorem3}
Let  $a\in C(\Gm)$ when $\Gm$ is a finite curve and $a\in
C(\dot{\Gamma})$ when  $\Gm$ is an infinite curve starting and
ending at infinity, where $\dot{\Gm}$ is the compactification of \
$\Gm$ by a single infinite point, that is,
$a(t(-\infty))=a(t(+\infty))$. Under the conditions of  Theorem
\ref{theorem2} the operator
$$ (S_\Gm a I - aS_\Gm)f= \frac{1}{\pi i}\intl_\Gm
\frac{a(\tau)-a(t)}{\tau-t}f(\tau)d\nu(\tau) $$  is compact in the
space $L^{p(\cdot)}(\Gm,\varrho)$.
\end{theorem}

In the Euclidean setting and  non-weighted case, the compactness
of the commutators generated by Calder\'on-Zygmund singular
operators in $\mathbb{R}^n$ with $b\in BMO$ ($b\in VMO$ in the
case of $\mathbb{R}^1$) was studied in \cite{299c} and
\cite{101zb}. The boundedness of multilinear commutators of
singular operators on variable exponent Lebesgue spaces was
studied in \cite{721z}.

 \vspace{3mm} The results obtained for the Cauchy
singular integral operator  led to a possibility to investigate
the Riemann problem in the setting of such spaces: find an
analytic function $\Phi$ on the complex plane cut along $\Gamma$
whose boundary values satisfy the conjugacy condition
\begin{equation}\label{f1}
\Phi^+(t)=G(t)\Phi^-(t)+g(t), \ \ \  t\in \Gm,
\end{equation}
where $G$ and $g$ are the given functions on $\Gm$ and $\Phi^+$
and $\Phi^-$ are boundary values of $\Phi$ on $\Gm$ from inside
and outside $\Gm$, respectively. The solution of (\ref{f1}) is
looked for in the class of analytic functions represented by the
Cauchy type integral with density in the spaces
$L^{p(\cdot)}(\Gm)$  and it is  assumed
 that $g$ belongs to the same class. We refer to \cite{171} and \cite{411z} for the classical solution
 of this problem and to \cite{63}, \cite{207}, \cite{208}, \cite{313}, \cite{629} for $L^p$-solutions with constant $p$.
 The investigation of the effects generated by the variable exponent setting of the  Riemann
 problem was given in \cite{317ab}, both for the case
when the coefficient $G$ is is continuous or piecewise continuous
and also for the case of oscillating coefficient. The Fredholmness
properties of the singular integral operators related to the
Riemann problem in the spaces $L^{p(\cdot)}(\Gm,\varrho)$ were
studied in \cite{321f} and in more general setting in
\cite{299ad}, \cite{299ae}.

\subsection{On potential operators }\label{subssobolev}
For non-weighted results on potentials and Sobolev embeddings  we
refer to \cite{101zb}, \cite{105a}, \cite{145zz}, \cite{147a},
\cite{160z}, \cite{167c}, \cite{405b}, \cite{577}.

 \vs{4mm} \textbf{a) Weighted
$p(\cdot)\to q(\cdot)$-boundedness; the Euclidean case.}

 The  known  generalization
 of Sobolev theorem by   Stein-Weiss for the case of power weights was extended in \cite{584a}, \cite{589b} to the variable exponent
 setting as follows.

\begin{theorem}\label{samsharvak}
Let $p\in \mathbb{P}(\rn)$,
$\sup\limits_{x\in\rn}p(x)<\frac{n}{\al}$,
$\varrho(x)=|x|^{\gm_0}(1+|x|)^{\gm_\infty-\gm_0}$ and
\begin{equation}\label{f:25ce1}
 |p_\ast(x)-p_\ast(y)|\le
\frac{A_\infty}{\ln \frac{1}{|x-y|}}, \ \
  \ \ \   |x-y|\leq \frac{1}{2}, \;\; x,\,y\in \mathbb{R}^n, \ \ \
  \ p_\ast(x)=p\left(\frac{x}{|x|^2}\right),
\end{equation}
 the operator
\begin{equation} \label{Riesz}
I^{\al} f(x) \ = \  \int\limits_{\rn}\frac{f(y)\;
dy}{|x-y|^{n-\al}} , \ \ \ \ \ \ 0<\al< n
\end{equation}
  is bounded from the space
$L^{p(\cdot)}(\mathbb{R}^n,\varrho)$ into the space
$L^{q(\cdot)}(\mathbb{R}^n,\varrho)$,  if
\begin{equation}\label{f7}
\al -\frac{n}{p(0)}< \gm_0<\frac{n}{p^\prime(0)}, \ \ \ \ \al -
\frac{n}{p(\infty)}< \gm_\infty<\frac{n}{p^\prime(\infty)}.
\end{equation}
\end{theorem}

\vs{4mm}Recently results similar to Theorem \ref{samsharvak} were
obtained for  more general weights like in Theorem \ref{theorem}.

We give  separate formulations of this generalizations for bounded
domains and for the whole space $\rn$ because for bounded domains
we are also able to admit variable orders $\al$ of the potential
operator.

We assume that
\begin{equation}\label{fvnuis4poniuska}
0 <\inf\limits_{x\in\Om}\al(x)p(x) \le
\sup\limits_{x\in\Om}\al(x)p(x)<n
\end{equation}
and
\begin{equation} \label{2.4abvNcbhriu}
|\al(x)-\al(y)|\le \frac{A}{\ln\frac{1}{|x-y|}} , \ \ \ \  x,y\in
\Om, \ \ \ \ \ \
 \ \ |x-y|\le \frac{1}{2}.
\end{equation}
\begin{theorem} \label{th:A}
Let $\Om$ be a bounded open set in $\rn$ and $x_0\in \ol{\Om}$,
let $p\in \mathbb{P}(\Om)$  and $\al$ satisfy conditions
(\ref{fvnuis4poniuska})-(\ref{2.4abvNcbhriu}). Let also
$\varrho(x)=w(|x-x_0|), x_0\in\ol{\Om}$, where
\begin{equation} \label{3.4Naga}
w(r)\in \Phi^\bt_\gm([0,\ell]) \ \ \ \   \textrm{with}   \ \ \ \
\bt=\al(x_0)-\frac{n}{p(x_0)}, \ \ \ \gm=\frac{n}{p^\prime(x_0)},
\end{equation}
or equivalently
\begin{equation} \label{3.5N}
w\in \widetilde{W}_0 \ \ \ \  \ \textrm{and} \ \ \ \ \
\al(x_0)-\frac{n}{p(x_0)} <m(w)\le M(w)<\frac{n}{p^\prime(x_0)}.
\end{equation}
 Then
\begin{equation} \label{3.3N}
\left\|I^{\al(\cdot)}f\right\|_{L^{q(\cdot)}\left(\Om,\varrho\right)}
\le C \left\|f\right\|_{L^{p(\cdot)}\left(\Om,\varrho\right)}.
\end{equation}
\end{theorem}

$w_0(r)$ belongs to some $\Phi_\gm^\bt$-class on $[0,1]$ and
 $w_\infty(r)$ belongs to some $\Psi_\gm^\bt$-class on $[1,\infty]$ and both
  the weights are continued by constant to $[0,\infty)$:
$$w_0(r)\equiv w_0(1), \ \ 1\le r<\infty  \ \ \  \textrm{and} \ \  \ w_\infty(r) \equiv w_\infty(1), \ \ \ 0< r\le 1 .$$

\begin{theorem} Let $0<\al<n$, $p\in
\mathbb{P}(\rn)$ satisfy
 assumption (\ref{f:25ce1}) and condition $p^+(\rn)<\frac{n}{\al}$, and let $\varrho(x)=
 w_0(|x|)w_\infty\left(|x|\right).$ The  operator $I^\al$ is bounded
from  $L^{p(\cdot)}(\mathbb{R}^n,\varrho)$ to
$L^{q(\cdot)}(\mathbb{R}^n,\varrho)$, if
\begin{equation} \label{3.4Naganew} w_0(r)\in \Phi^{\bt_0}_{\gm_0}([0,1]), \ \  w_\infty(r)\in
\Psi^{\bt_\infty}_{\gm_\infty}([1,\infty))
\end{equation}
where $ \bt_0=\al -\frac{n}{p(0)}, \ \gm_0=\frac{n}{p^\prime(0)},
\ \ \bt_\infty=\frac{n}{p^\prime(\infty)}, \ \gm_\infty=\al
-\frac{n}{p(\infty)},
 $
or equivalently
\begin{equation} \label{3.5Nnew}
w_0\in \widetilde{W}([0,1]),  \ \ \al -\frac{n}{p(0)} <m(w_0)\le
M(w_0)< \frac{n}{p^\prime(0)},
\end{equation}
and
\begin{equation} \label{3.5Nnew1}
w_\infty \in \widetilde{W}([1,\infty]),  \ \ \ \al
-\frac{n}{p(\infty)}<m(w_\infty)\le
M(w_\infty)<\frac{n}{p^\prime(\infty)}.
\end{equation}
\end{theorem}
We refer also to the paper \cite{145c} where there was obtained an
extension of the Adams' trace inequality for potentials to the
variable exponent setting on homogeneous spaces.

\vs{4mm} \textbf{b) Weighted $p(\cdot)\to q(\cdot)$-boundedness
for potentials on Carleson curves.}
 \vs{4mm} A statement on $p(\cdot)\to q(\cdot)$
boundedness of potential operators
\begin{equation} \label{potential}
I^{\al(\cdot)} f(t) \ = \  \int\limits_{\Gm}\frac{f(\tau)\;
d\nu(\tau)}{|t-\tau|^{1-\al(t)}}
\end{equation}
 is also valid on an arbitrary Carleson curve $\Gm$
and for variable orders $\al(t)$ as well, see \cite{321j}, as
given below.

\begin{theorem}\label{theorem4} Let
 $\Gm$ be a simple Carleson curve of a finite  length,
 $p\in \mathbb{P}(\Gm)$ and  $\al(t)$ satisfy the assumptions
\begin{equation}\label{fvnuis4}
0 <\inf\limits_{t\in\Gm}\al(t)p(t) \le
\sup\limits_{t\in\Gm}\al(t)p(t)<1.
\end{equation}
Then the operator $I^{\al(\cdot)}$ is bounded from the space
$L^{p(\cdot)}(\Gm)$ into the space $L^{q(\cdot)}(\Gm)$ with
$\frac{1}{q(t)}=\frac{1}{p(t)}-\al(t)$. This statement remains
valid for infinite Carleson curves if, in addition to the above
conditions, $p(t)=p_\infty=const$ and $\al(t)=\al_\infty=const$
outside some circle $B(t_0,R), t_0\in \Gm$.
\end{theorem}

\vspace{2mm} In the next weighted version  of Theorem
\ref{theorem4} for finite curves we additionally suppose that the
order $\al(t)$ is log-continuous at the nodes of the weight:
\begin{equation}\label{poch}
|\al(t)-\al(t_k)|\le \frac{A}{|\ln|t-t_k||}, \quad k=1,...,m.
\end{equation}

\begin{theorem}\label{theorem5}
 Let
 $\Gm$ be a simple Carleson curve of a finite  length.
 Under  condition (\ref{poch}) and the assumptions  of Theorem
 \ref{theorem4},
 the operator $I^{\al(\cdot)}$ is bounded from the space
$L^{p(\cdot)}(\Gm,\varrho)$ into the space
$L^{q(\cdot)}(\Gm,\varrho)$ where \
$\frac{1}{q(t)}=\frac{1}{p(t)}-\al(t)$, and  the weight
$\varrho(t) = \prod\limits_{k=1}^m |t-t_k|^\bt_k $, if
\begin{equation}\label{f34a}
\al(t_k)-\frac{1}{p(t_k)}<\bt_k< 1-\frac{1}{p(t_k)}, \ \ \ \ \
k=1,...,m.
\end{equation}
\end{theorem}

\begin{corollary}  Under the assumptions of Theorem \ref{theorem5}, the
fractional maximal operator
$$M_{\al(\cdot)}f(t)= \sup\limits_{r>0}\frac{1}{\nu\{\Gm(t,r)\}^{n-\al(t)}}\int\limits_{\Gm(t,r)}|f(\tau)|\;
d\nu(\tau), \quad \quad \Gm(t,r)= \Gm \cap B(t,r)
$$
is bounded from the space $L^{p(\cdot)}(\Gm,\varrho)$ into the
space $L^{q(\cdot)}(\Gm,\varrho)$.
\end{corollary}

We mention also Hardy-type inequalities for potentials obtained in
\cite{580b}, \cite{580bc} in the multidimensional case and in
\cite{146za} in the one-dimensional case,  including in particular
the case of Riemann-Liouville and Weyl fractional integrals.

\vs{4mm} \textbf{c) Characterization of the range of potential
operators.}
 The inversion of the Riesz potentials with densities
in $L^{p(\cdot)}(\rn)$ by means of hypersingular integrals
$\mathbb{D}^\al f$ (Riesz fractional derivatives of order $\al$)
was obtained in \cite{18b} (we refer to \cite{580} for the case of
constant $p$ and hypersingular integrals in general).  This gave a
possibility to give   in \cite{18d} a characterization of the
range $I^\al[L^{p(\cdot)}(\rn)]$
 in terms of convergence  of $\mathbb{D}^\al f$ in
 $L^{p(\cdot)}(\rn)$ as follows.

\begin{theorem} \label{rieszrange1}
Let $p\in \mathbb{P}(\rn)$, $0< \alpha <n$, $1< p_-(\rn) \leq
p^+(\rn) < \frac{n}{\alpha}$ and let $f$ be a locally integrable
function. Then $f\in I^\alpha [ L^{p(\cdot)}(\rn) ]$, if and only
if $f\in L^{q(\cdot)}(\rn)$, with $\frac{1}{q(\cdot)} =
\frac{1}{p(\cdot)} - \frac{\alpha}{n}$, and there exists the Riesz
derivative $\mathbb{D}^\alpha f$ (in the sense of convergence in
$L_{p(\cdot)}$).
\end{theorem}

A study of the range $I^\alpha [ L^{p(\cdot)}(\Om) ]$ for domains
$\Om\subset \rn$ is an open question; in the form  given in
Theorem \ref{rieszrange1} it is open even in the case of constant
$p$, one of the reasons being in the absence of the corresponding
apparatus of hypersingular integrals adjusted  to domains in
$\rn$; some their analogue reflecting the influence of the
boundary was recently suggested in \cite{504z} for the case
$0<\al<1.$  In the one dimensional case for $\Om=(a,b),
-\infty<a<b\le\infty$, when the range of the potential coincides
with that of the Riemann-Liouville fractional integral operators
(in the case $1<p^+(a,b)<\frac{1}{\al}$), the characterization for
variable $p(x)$ was obtained in \cite{504za}, where for
$-\infty<a<b<\infty$ there was also shown its coincidence with the
space of restrictions of Bessel potentials.

  The result of Theorem \ref{rieszrange1} was used in \cite{18d}
  to obtain a characterization  of the
\emph{Bessel potential space}
$$
\mathcal{B}^\alpha [L^{p(\cdot)}(\rn)] = \{ f: \; f=
\mathcal{B}^\alpha \varphi, \quad \varphi \in
L^{p(\cdot)}(\mathbb{R}^n) \}, \quad \alpha \geq 0,
$$
where $\mathcal{B}^\al\vi = F^{-1}(1+|\xi|^2)^{-\al/2}F\vi$ and
runs as follows.

\begin{theorem} \label{Bessel}
Under the conditions of Theorem \ref{rieszrange1}
\begin{equation}\label{formula}
\mathcal{B}^\alpha [L^{p(\cdot)}(\rn)] = L^{p(\cdot)}(\rn)\bigcap
I^\alpha [ L^{p(\cdot)} (\rn)] = \{f\in L^{p(\cdot)}(\rn): \
\mathbb{D}^\al f \in L^{p(\cdot)}(\rn)\}
\end{equation}
 and
$\mathcal{B}^m [L^{p(\cdot)}(\rn)] = W^{m,p(\cdot)}(\rn)$ for any
integer $m\in \mathbb{N}_0$, where $W^{m,p(\cdot)}(\rn)$ is the
Sobolev space with the variable exponent $p(x)$.
\end{theorem}

\vs{4mm} Statement (\ref{formula}) has the following
generalization, see \cite{504za}, Theorem 4.10. (We refer to
\cite{53b} for the notion of Banach function spaces).

\begin{theorem}\label{theoalmbessel}
Let $\mathrm{Y}=\mathrm{Y}(\mathbb{R}^n)$ be a Banach function
space, satisfying the assumptions \\ i)\ $C^\infty_{0}$ is dense
in $\mathrm{Y}$; \\ ii) \  the maximal operator $\mathcal{M}$ is
bounded in $\mathrm{Y}$; \\ iii) \ $I^\al f(x)$ converges
absolutely for
almost all $x$ for every $f\in\mathrm{Y}$ and \\
$(1+|x|)^{-n-\al}I^\al f(x) \in L^1(\mathbb{R}^n)$.

\noindent Then
\begin{equation}\label{almbessel}
\mathcal{B}^\al(\mathrm{Y})=\mathrm{Y} \bigcap I^\al(\mathrm{Y})=
\{f\in \mathrm{Y}: \mathbb{D}^\al f=\lim\limits_{\varepsilon
\rightarrow 0\atop (\mathrm{Y})}\mathbb{D}_\ve^\al f\in \mathrm{Y}
\}.
\end{equation}
\end{theorem}

From Theorem \ref{theoalmbessel} there follows, in particular, the
characterization of the ranges of potential operators over
weighted Lebesgue spaces with variable exponent obtained by means
of results of Subsection \ref{subsmaxop}  for the maximal
operator.

Observe that certain results related to imbedding of the range of
the Riesz potential operator into  H\"older spaces (of variable
order) in the case $p(x)\ge n$ were obtained in  \cite{18e}. The
results proved in \cite{18e} run as follows where
\begin{equation} \label{pi}
\Pi_{p,\Omega}:=\{x\in\Omega: p(x)>n\}
\end{equation}
and $C^{0,\alpha(\cdot)}(\Omega)$ is  the space of bounded
continuous  functions $f$ with a finite seminorm
$$
[f]_{\alpha(\cdot),\Omega} := \sup\limits_{x, x+h \in \Omega \atop
0<|h|\leq 1} \frac{|f(x+h)-f(x)|}{|h|^{\alpha(x)}}.
$$

\begin{theorem} \label{pointwise}
Let $\Omega$ be a bounded open set with Lipschitz boundary and let
 $p(\cdot)$ satisfy the  log-condition (\ref{2nd})
and have a non-empty set $\Pi_{p,\Omega}$. If $f\in
W^{1,p(\cdot)}(\Omega)$, then
\begin{equation} \label{pointwise1}
|f(x)-f(y)| \leq C(x,y)\; \||\nabla f|\|_{p(\cdot),\Omega} \;
|x-y|^{1-\frac{n}{\min[p(x),p(y)]}}
\end{equation}
for all $x,y\in \Pi_{p,\Omega}$ such that $|x-y|\le 1$, where
$$C(x,y)=\frac{c}{\min[p(x),p(y)]-n} $$ with $c>0$ not depending on
$f,x$ and $y$.
\end{theorem}

\begin{theorem} \label{embeddomain}
Let $\Omega$ be a bounded open set with Lipschitz boundary and
suppose that $p(\cdot)$ satisfies the logarithmic condition
(\ref{2nd}). If $\inf\limits_{x\in\Omega}p(x)>n$, then
\begin{equation} \label{embedvarorder1}
W^{1,p(\cdot)}(\Omega) \hookrightarrow
C^{0,1-\frac{n}{p(\cdot)}}(\Omega),
\end{equation}
where ``$\hookrightarrow$" means continuous embedding.
\end{theorem}

Theorem \ref{embeddomain} is an improved version of the result
earlier obtained in \cite{147a}, \cite{163a}. The papers
\cite{167ac}, \cite{167ab}  are also relevant to the topic. We
refer also to \cite{224z} where the capacity approach was used to
get embeddings into the space of continuous functions or into
$L^\infty(\Omega)$.

 Observe that is \cite{18e} there
were also obtained $W^{1,p(\cdot)}(\Om)\to L^{q(\cdot)}$-estimates
of hypersingular integrals (fractional differentiation operators)
\begin{equation} \label{defhypersing}
\mathcal{D}^{\alpha(\cdot)} f(x) = \int_{\Omega}
\frac{f(x)-f(y)}{|x-y|^{n+\alpha(x)}} \, dy, \quad x\in \Omega.
\end{equation}

 \vs{4mm}

\subsection{On Hardy operators }\label{Hardy}

For Hardy operators (\ref{Hardy}) in \cite{107da} the following
result was obtained without the log condition at all points on
$\mathbb{R}^1_+$.  By $\mathcal{M}_{0,\infty}(\mathbb{R}^1_+)$
 we denote the set  of all measurable bounded functions $p(x): \mathbb{R}^1_+\to \mathbb{R}^1_+$
which satisfy the following conditions:

\vspace{2mm}\noindent
 $ i) \ \ \ \ \ \ 0\le p_-\le
p(x)\le p_+<\infty, \ \ \ \ x\in
\mathbb{R}^1_+$,\\
$ ii_0) \ \ $  there exists  $p(0)=\lim\limits_{x\to 0}p(x)$ and \
$|p(x)-p(0)|\le \frac{A}{\ln\,\frac{1}{x}}, \ \ 0<x\le
\frac{1}{2},$\\
 $ ii_\infty)$ \ there exists $\mu(\infty)=\lim\limits_{x\to \infty}p(x)$ \ \ \
\textrm{and} \ $|p(x)-p(\infty)|\le \frac{A}{\ln\,x},  \  \ x\ge
2. $

By $\mathcal{P}_{0,\infty}=\mathcal{P}_{0,\infty}(\mathbb{R}^1_+)$
we denote the subset of functions $p(x)\in
\mathcal{M}_{0,\infty}(\mathbb{R}^1_+)$ with
$\inf\limits_{x\in\mathbb{R}^1_+}p(x)\ge 1$.
\begin{theorem}\label{7778}  Let $p,q\in \mathcal{P}_{0,\infty}(\mathbb{R}^1_+)$
and $\mu\in\mathcal{M}_{0,\infty}(\mathbb{R}^1_+)$  and let
 $$
 \frac{1}{q(0)}=\frac{1}{p(0)}-\mu(0), \  \frac{1}{q(\infty)}=\frac{1}{p(\infty)}-\mu(\infty) \quad
 \textrm{and}\quad  0\le \mu(0)
< \frac{1}{p(0)}, \  0\le \mu(\infty) < \frac{1}{p(\infty)}.$$
 Then the Hardy-type inequalities
\begin{equation}\label{citx1new}
\left\|x^{\al+\mu(x)-1}\intl_0^x\frac{f(y)\,dy}{y^\al}\right\|_{L^{q(\cdot)}(\mathbb{R}^1_+)}
\le C \left\|f\right\|_{L^{p(\cdot)}(\mathbb{R}^1_+)}
\end{equation}
and
\begin{equation}\label{citdox2new}
\left\|x^{\bt+\mu(x)}\intl_x^\infty\frac{f(y)\,dy}{y^{\bt
+1}}\right\|_{L^{q(\cdot)}(\mathbb{R}^1_+)} \le C
\left\|f\right\|_{L^{p(\cdot)}(\mathbb{R}^1_+)},
\end{equation}
are valid, if and only if
\begin{equation} \label{0v08cv}
\al <  \min\left\{\frac{1}{p^\prime(0)},
\frac{1}{p^\prime(\infty)}\right\} \quad \textrm{and}\quad \bt >
\max\left\{\frac{1}{p(0)},\frac{1}{p(\infty)}\right\}.
\end{equation}
\end{theorem}

For previous version of Hardy inequality we refer to \cite{224zc},
\cite{ 321a}. In \cite{224zc} a  multidimensional  analogue of
Hardy inequality was also considered.

\section{ Weighted boundedness of maximal operators on metric measure spaces}\label{yyyys4}
\setcounter{theorem}{0} \setcounter{equation}{0}

 In
the case of constant $p\in (1,\infty)$ the boundedness of  the
maximal operator on bounded metric measure  spaces is well known,
due to A.P.Calder\'on  \cite{72b} and R.Mac\'{i}as and C.Segovia
\cite{381a} for weights in  the Muckenhoupt class
$\mathcal{A}_p=\mathcal{A}_p(X)$, defined by the condition
\begin{equation}\label{ap}
\sup\limits_{x\in X, r>0}\left(\frac{1}{\mu
B(x,r)}\intl_{B(x,r)}|\varrho(y)|^{p}d\mu(y)\right)\left(\frac{1}{\mu
B(x,r)}\intl_{B(x,r)}
\frac{d\mu(y)}{|\varrho(y)|^{p^\prime}}\right)^{p-1} <\infty.
\end{equation}

For variable exponents the maximal operator on metric measure
spaces  was considered
 in \cite{224b} and \cite{307b}, where the following non-weighted result was obtained.

\begin{theorem}\label{khabazi-hasto}
Let a bounded metric measure  space  $X$ satisfy the doubling
condition (\ref{double}) and  $p\in \mathbb{P}(X)$. Then the
maximal operator $\mathcal{M}$ is bounded in the space
$L^{p(\cdot)}(X)$.
\end{theorem}

As was observed in \cite{224b}, in contrast to the case of
constant $p$, the doubling condition is not necessary for the
boundedness of the maximal operator when $p$ is variable.

 \vs{4mm} The boundedness of the operator
$\mathcal{M}$ in the weighted space  $L^{p(\cdot)}(X,\varrho)$ is
known for the cases where $X=\Om$ is a bounded domain in
$\mathbb{R}^n$  or $X=\Gm$ is a Carleson curve on the complex
plane,  see Theorems \ref{theorem} and \ref{theorem1}.

\vs{4mm}
 For an arbitrary  metric measure space with doubling condition,
 we present in this section new results on  weighted boundedness
  given in Theorems A, B and C stated below. There proof taking
  too much space will be given elsewhere.

Let $\mathcal{A}_{p(\cdot)}(X)$ be the class (\ref{muckenh}). To
formulate Theorem A we introduce the following "Muckenhoupt-like
looking" class
 $\widetilde{\mathcal{A}}_{p(\cdot)}(X)$ of weights, which
satisfy the condition
\begin{equation}\label{erzatz}
\sup\limits_{x\in X, r>0}\left(\frac{1}{\mu
B(x,r)}\intl_{B(x,r)}|\varrho(y)|^{p(y)}d\mu(y)\right)\left(\frac{1}{\mu
B(x,r)}\intl_{B(x,r)}
\frac{d\mu(y)}{|\varrho(y)|^\frac{p(y)}{p_--1}}\right)^{p_--1}
<\infty.
\end{equation}
This class $\widetilde{\mathcal{A}}_{p(\cdot)}(X)$ used in Theorem
A is narrower than the class $\mathcal{A}_{p(\cdot)}$. However, it
coincides with the Muckenhoupt class $\mathcal{A}_p$ in case $p$
is constant.

In Theorem A, under log-condition on $p$ and doubling condition on
the measure we show that
$$\widetilde{\mathcal{A}}_{p(\cdot)}(X)\subset {\mathcal{A}}_{p(\cdot)}(X).$$
 In Theorem B we deal with a special class of radial type weights
in the Zygmund-Bary-Stechkin class and arrive at the necessity to
relate the properties of the weight to those of the measure $\mu
B(x,r)$ as stated in (\ref{avotyavamhu}). Such a result for the
Euclidean case was earlier obtained in \cite{317c}. The proof for
the case of metric measure spaces requires an essential
modification of the technique used. Theorem B is proved by means
of Theorem A, but it is not contained in Theorem A, being more
general in its range of applicability.

\vspace{3mm}

 \vspace{4mm}\textbf{Theorem A.} \
\textit{Let $X$ be a bounded doubling metric measure space, let
the exponent $p\in \mathbb{P}(X)$ and the weight $\varrho$ fulfill
condition (\ref{erzatz}). Then the operator $\mathcal{M}$ is
bounded in $L^{p(\cdot)}(X,\varrho)$.}

\vspace{3mm}

In Theorems B  and C we deal with bounded and unbounded metric
spaces, respectively. In Theorem B we consider weights of the form
\begin{equation}\label{u56741}
\varrho(x)=\prod_{k=1}^N w_k(d(x,x_k)), \ \ \ x_k\in X,
\end{equation}
 where  $x_k$ are distinct points and  $w_k(r)$
 may oscillate between two power functions  as $r\to 0+$
(radial Zygmund-Bary-Stechkin type weights), and in Theorem C  we
consider  similar weights of the form
\begin{equation}\label{fdbs6}
\varrho(x)= w_0[1+d(x_0,x)]\prod\limits_{k=1}^N w_k[d(x,x_k)], \ \
\ x_k\in X, k=0,1,...N.
\end{equation}

We make also use of the following   numbers which play a role of
dimensions of the space $(X,d,\mu)$ at the
point $x\in X$:\\
 1) the \textit{ local lower and upper dimensions}
\begin{equation}\label{i}
m(\mu B_x)= \sup_{t>1} \frac{\ln\left(\liminf\limits_{r\to 0}
\frac{\mu B(x,rt)}{\mu B(x,r)}\right)}{\ln t},  \quad M(\mu B_x)=
\inf_{t>1} \frac{\ln\left(\limsup\limits_{r\to 0} \frac{\mu
B(x,rt)}{\mu B(x,r)}\right)}{\ln t},
\end{equation}
2)  \textit{similar dimensions "influenced" by infinity}:
\begin{equation}\label{ia}
m_\infty(\mu B)= \sup_{t>1} \frac{\ln\left(\liminf\limits_{r\to
\infty} \frac{\mu B(x,rt)}{\mu B(x,r)}\right)}{\ln t}, \quad
M_\infty(\mu B)= \inf_{t>1} \frac{\ln\left(\limsup\limits_{r\to
\infty} \frac{\mu B(x,rt)}{\mu B(x,r)}\right)}{\ln t}
\end{equation}
the latter appearing only in the case of unbounded $X$. The idea to use the above local dimensions
was borrowed from  papers \cite{539j}, \cite{539jnew}. It may be shown that the numbers
$m_\infty(\mu B)$ and $M_\infty(\mu B)$ do not depend on $x$, see \cite{539j}, \cite{539jnew}.
\begin{remark}\label{remark}
 In a different form  local dimensions were introduced
and/or used in \cite{160zzzz}, \cite{160zzzza}, \cite{224ab}, \cite{224b}, \cite{244a}. The
introduction of the local dimensions in the form described above was influenced  by the study of
lower and upper indices of oscillating almost increasing functions in \cite{539}, \cite{539e},
\cite{539d} and  application of that study in  \cite{317d}, \cite{317c}, \cite{317e}, \cite{539j},
\cite{539jnew}, \cite{539h}.
\end{remark}

It may be shown (see \cite{539h}) that  for an arbitrarily small
$\ve>0$
\begin{equation}\label{iihochejk}
c_1 r^{M(\mu B_x)+\ve}\le \mu B(x,r) \le  c_1 r^{m(\mu B_x)-\ve},
\quad 0< r \le R<\infty
\end{equation}
and
\begin{equation}\label{i0hejk}
c_3 r^{m_\infty(\mu B)-\ve}\le \mu B(x,r) \le  c_4 r^{M_\infty(\mu
B)+\ve}, \quad  r_0\le r < \infty,
\end{equation}
where $c_i,i=1,2,3,4$, depend on $\ve>0$, but do not depend on $r$
and $x$.

 \vspace{3mm} In the sequel we will use the "uniform" index $m(\mu B)$ introduced in
\eqref{introduced}.

The Zygmund-Bary-Stechkin class $\Phi_1^0$ of weights and the
upper and lower indices of weights (of the type of
Matuszewska-Orlicz indices, see \cite{382a}, close in a sense to
the Boyd indices) used in the theorem below, were defined in
Section \ref{Phi}. Various non-trivial examples of functions in
Zygmund-Bary-Stechkin-type classes with coinciding indices may be
found in \cite{539}, Section II; \cite{539a}, Section 2.1, and
with non-coinciding indices in \cite{539d}.

\vskip+0.2cm \textbf{Theorem B}. \  {\em Let $X$ be a bounded
doubling metric measure  space and let $p\in \mathbb{P}(X)$. The
operator $\mathcal{M}$ is bounded in $L^{p(\cdot)}(X,\varrho)$
with weight (\ref{u56741}), if $r^{\frac{m(\mu
B)}{p(x_k)}}w_k(r)\in \Phi_{m(\mu B)}^0$}, or equivalently $w_k
\in \widetilde{W}([0,\ell]),  \  \ell=\diam X,$ and
\begin{equation}\label{ere}
-\frac{m(\mu B)}{p(x_k)} < m(w_k)\le M(w_k) < \frac{m(\mu
B)}{p^\prime(x_k)}\, \ , \ \ k=1,2,...,N.
\end{equation}

In the case where $X$ is a bounded open set in $\mathbb{R}^n$,
Theorem B was proved  in \cite{317c} and coincides with Theorem
\ref{theorem}; in the case where $X=\Gm$ is a Carleson curve, it
was proved in \cite{321j} for  power weights, as stated in Theorem
\ref{theorem1}, and in \cite{317d} for  weights in
Zygmund-Bary-Stechkin type class.

\vspace{3mm} \textbf{Theorem C}. \ \textit{Let $X$ be an unbounded
doubling metric measure  space and let $p\in \mathbb{P}(X)$, and
let there exist a ball $B(x_0,R), x_0\in X$ such that $p(x)\equiv
p_\infty =const$ for $x\in X\backslash B(x_0,R)$. Then the maximal
operator $\mathcal{M}$ is bounded in the space
$L^{p(\cdot)}(X,\varrho)$, with weight (\ref{fdbs6}), if $w_k \in
\widetilde{W}(\mathbb{R}_+^1)$ and
\begin{equation}\label{f27d}
-\frac{m({\mu B})}{p(x_k)}<m(w_k)\le M(w_k)<\frac{m(\mu
B)}{p^\prime(x_k)}, \ k=1,...,N,
\end{equation}
and
\begin{equation}\label{f27dco}
 -\frac{m_\infty(\mu
B)}{p_\infty}<\sum\limits_{k=0}^N m_\infty(w_k)\le
\sum\limits_{k=0}^N M_\infty(w_k) <\frac{m_\infty(\mu
B)}{p^\prime_\infty}-\Delta_{p_\infty},
\end{equation}
where $\Delta_{p_\infty}=\frac{M_\infty(\mu B)-m_\infty(\mu
B)}{p_\infty}.$ }

\vspace{3mm}In particular, for the power type weight
\begin{equation}\label{fdbs6co}
\varrho(x)= (1+d(x_0,x))^{\bt_0}\prod\limits_{k=1}^N
[d(x,x_k)]^{\bt_k}, \ \ \  x_k\in X, k=0,1,...N
\end{equation}
conditions (\ref{f27d})-(\ref{f27dco}) take the form
\begin{equation}\label{f27dcopra}
-\frac{m({\mu B})}{p(x_k)}<\bt_k<\frac{m(\mu B)}{p^\prime(x_k)}, \
\quad k=1,...,N,
\end{equation}
and
\begin{equation}\label{f27dcoprasiuk}
 -\frac{m_\infty(\mu
B)}{p_\infty}<\sum\limits_{k=0}^N \bt_k <\frac{m_\infty(\mu
B)}{p^\prime_\infty}-\Delta_{p_\infty}.
\end{equation}

The bounds  in (\ref{f27dco}) turn to take a natural form
$-\frac{m_\infty(\mu B)}{p_\infty}$ and $\frac{m_\infty(\mu
B)}{p^\prime_\infty}$ with $\Delta_{p_\infty}=0$ when "dimensions"
$m_\infty(\mu B)$ and $M_\infty(\mu B)$ coincide with each other.
In particular, in the case where $X$ has a constant dimension
$d>0$ in the sense that
$$C_1r^d\le \mu B(x,r)\le C_2 r^d,$$
conditions (\ref{f27dcopra})-(\ref{f27dcoprasiuk}) take the form
\begin{equation}\label{f27dcopra1n}
-\frac{d}{p(x_k)}<\bt_k<\frac{d}{p^\prime(x_k)}, \ \quad
k=1,...,N, \quad  -\frac{d}{p_\infty}<\sum\limits_{k=0}^N \bt_k
<\frac{d}{p^\prime_\infty}.
\end{equation}

 \vspace{2mm} The Euclidean space version of Theorem  C for variable exponents and power weights was obtained in
  \cite{307a}.

It goes without saying that in Theorems B and C we should be interested in the result with  the
values $m({\mu  B_{x_k}})$ instead of taking a kind of the  infinum in $x$ of $m({\mu  B_x})$ with
respect to all $x$. However, such a localization of the index  $m(\mu B_x)$ remains an open
question. So, instead, we deal with that infimum which, of course, serves for all the points $x_k$.

\section*{Acknowledgments}

 This work was made under the project "Variable Exponent Analysis" supported by INTAS grant
 Nr.06-1000017-8792.

%\bibliographystyle{plain}

%\bibliography{bibl_198}
\def\ocirc#1{\ifmmode\setbox0=\hbox{$#1$}\dimen0=\ht0 \advance\dimen0
  by1pt\rlap{\hbox to\wd0{\hss\raise\dimen0
  \hbox{\hskip.2em$\scriptscriptstyle\circ$}\hss}}#1\else {\accent"17 #1}\fi}

\end{document}